\documentclass[a4paper,UKenglish,cleveref, autoref, thm-restate]{lipics-v2021}
\nolinenumbers

\usepackage{hyperref}

\usepackage{amsmath}
\usepackage{upgreek}
\usepackage{xcolor}
\usepackage{colortbl}
\usepackage[capitalise]{cleveref}

\usepackage{pbox}
\usepackage{makecell}

\usepackage{tikz-cd}
\usepackage{tabularx}
\definecolor{keywordcolor}{rgb}{0.7, 0.1, 0.1}   
\definecolor{tacticcolor}{rgb}{0.0, 0.1, 0.3}    
\definecolor{commentcolor}{rgb}{0.4, 0.4, 0.4}   
\definecolor{stringcolor}{rgb}{0.5, 0.3, 0.2}    
\definecolor{symbolcolor}{rgb}{0.1, 0.2, 0.7}    
\definecolor{sortcolor}{rgb}{0.1, 0.5, 0.1}      
\definecolor{attributecolor}{rgb}{0.7, 0.1, 0.1} 
\definecolor{errorcolor}{rgb}{1, 0, 0}           

\usepackage{listings}
\usepackage{flushend}
\usepackage{xspace}
\usepackage{underscore}
\usepackage{microtype}
\usepackage{datetime}
\usepackage{MnSymbol} 
\usepackage{ragged2e}
\newtheorem{assumption}{Assumption}
\newcommand{\lean}[1]{\lstinline[language=lean]{#1}}
\newcommand{\mathlib}{\textsf{mathlib}\xspace}

\usepackage{scalefnt}
\title{Formalization of Algorithms for Optimization with Block Structures} 

\author{Chenyi Li}{School of Mathematical Sciences, Peking University, China}{lichenyi@stu.pku.edu.cn}{}{}
\author{Zichen Wang}{School of Mathematics and Statistics, Xi’an Jiaotong University, China}{princhernwang@gmail.com}{}{}
\author{Yifan Bai}{School of Mathematical Sciences, Shanghai Jiao Tong University, China}{yifanbai@sjtu.edu.cn}{}{}
\author{Yunxi Duan}{School of Mathematics, Statistics and Mechanics, Beijing University of Technology, China}{yunxiduan@emails.bjut.edu.cn}{}{}
\author{Yuqing Gao}{Institut für Numerische Simulation, Universität Bonn, Germany}{gao.yuqing@uni-bonn.de}{}{}
\author{Pengfei Hao}{School of Mathematical Sciences, Beijing Normal University, China}{202111998147@mail.bnu.edu.cn}{}{}
\author{Zaiwen Wen\footnote{corresponding author}}{Beijing International Center for Mathematical Research, Peking University, China}{wenzw@pku.edu.cn}{}{}

\authorrunning{C.\,Li, Z.\,Wang, Y.\,Bai, Y.\,Duan, Y.\,Gao, P.\,Hao, and Z.\,Wen}

\Copyright{Chenyi Li, Zichen Wang, Yifan Bai, Yunxi Duan, Pengfei Hao, Yuqing Gao, and Zaiwen Wen} 



\category{} 

\relatedversion{} 

\supplement{Our code for this paper can be found at \url{https://github.com/optsuite/optlib}}

\acknowledgements{We want to thank the \mathlib community for useful discussion on this project. We also want to thank the students from the summer school and winter school at Peking University.}

\begin{document}

\maketitle

\begin{abstract}
Block-structured problems are central to advances in numerical optimization and machine learning. This paper provides the formalization of convergence analysis for two pivotal algorithms in such settings: the block coordinate descent (BCD) method and the alternating direction method of multipliers (ADMM). Utilizing the type-theory-based proof assistant Lean4, we develop a rigorous framework to formally represent these algorithms. Essential concepts in nonsmooth and nonconvex optimization are formalized, notably subdifferentials, which extend the classical differentiability to handle nonsmooth scenarios, and the Kurdyka–Łojasiewicz (KL) property, which provides essential tools to analyze convergence in nonconvex settings. Such definitions and properties are crucial for the corresponding convergence analyses. We formalize the convergence proofs of these algorithms, demonstrating that our definitions and structures are coherent and robust. These formalizations lay a basis for analyzing the convergence of more general optimization algorithms.

\end{abstract}

\section{Introduction}
Large-scale optimization problems with block structures play a central role in numerous fields, including machine learning, signal processing, and image processing. Efficiently solving these problems is essential for theoretical advancement and practical applications. Two prominent algorithms addressing these challenges are the block coordinate descent (BCD) method and the alternating direction method of multipliers (ADMM). The BCD algorithm solves optimization problems with objective functions comprising both separable and nonseparable components by updating the variable blocks by blocks. The ADMM algorithm primarily addresses optimization problems involving multi-blocks connected with coupling constraints, iteratively minimizing the augmented Lagrangian function for each block and updating the multipliers to enforce constraint satisfaction. Despite their methodological differences, both approaches share the fundamental strategy of decomposing complex optimization problems into simpler subproblems.

The convergence analysis of algorithms for optimization problems with block structures has been widely studied. For the BCD method, Powell \cite{powell1973search} demonstrated the challenges of nonconvex optimization through examples of non-convergence in nonconvex objective functions. Tseng \cite{tseng2001convergence} proposes a robust convergence framework for nondifferentiable convex minimization problems. To address these complexities in nonsmooth and nonconvex analysis, additional conditions, such as the Kurdyka-Łojasiewicz (KL) property, become necessary. The KL property can trace back to Łojasiewicz’s work on gradient inequalities \cite{lojasiewicz1963propriete}. The convergence of the BCD algorithm in nonconvex settings with the KL property is given by \cite{bolte2014proximal}. The ADMM method also gains significant attention due to its applicability in high-dimensional data problems and distributed optimization \cite{boyd2011distributed}. 
The classical ADMM method is introduced in \cite{glowinski1975approximation} and \cite{gabay1976dual}. Glowinski provides an introduction to the historical development in \cite{Glowinski2014}. The convergence analysis of the ADMM algorithm under certain conditions is established in \cite{gabay1976dual, glowinski2013numerical, fortin2000augmented}.  Fankel et al. give the convergence proof of a computationally more attractive large step-length version of the ADMM method in \cite{Fazel2013hankelmatrix}.

Several efforts have been made to formalize numerical algorithms using various formalization languages such as Coq \cite{huet1997coq}, Isabelle \cite{Nipkow2002APA}, and Lean \cite{de2015lean}. Tekriwal et al. formalize the convergence of the Jacobi method for linear systems in Coq \cite{Tekriwal2023Verified}, and the formalization of asymptotic convergence for stationary iterative methods such as Gauss-Seidel method is explored in \cite{Tekriwal2024Formalizatioin}. The analysis of floating-point errors in numerical algorithms is formalized in Coq \cite{Appel2024VCFloat, Kellison2024Numerical}. Additionally, the formalization of convergence rates for various first-order optimization methods is presented in Lean4 \cite{li2024formalization}. 

The Lean4 language and its theorem library \mathlib \cite{mathlibcommunity} offers a flexible formalization tool for mathematicians, enabling precise reasoning about mathematical concepts through a machine-verified framework. Lean4 is based on the dependent type theory, allowing for the expression of refined mathematical properties directly in the type system. In this setting, proofs are constructed as objects of a specific type, ensuring that every step is logically sound and verifiable. This approach enables the construction of rigorous formalizations across various fields of mathematics. For example, the Gagliardo-Nirenberg-Sobolev inequality is formalized by Doorn and Macbeth \cite{vandoorn_et_al:LIPIcs.ITP.2024.37}. Dupuis et al. explore the formalization of functional analysis with semilinear maps \cite{dupuis_et_al:LIPIcs.ITP.2022.10}. Moreover, a formalization of Doob’s martingale convergence theorems in mathlib is provided by Ying et al. \cite{Ying2023Doob}. 

We presents a formalized convergence analysis of the BCD and ADMM algorithms within the Lean4 theorem proving environment. Our contributions are structured as follows.
\begin{enumerate}
    \item \textbf{Subdifferential theory in Hilbert spaces.} The Fréchet and limiting subdifferentials are formalized, extending the formalization of the differential structure of nonsmooth functions. Basic theorems are presented, including equivalent definitions, subdifferential calculus, and optimality conditions of nonsmooth unconstrained optimization.
    \item \textbf{The KL property and KL functions.} The formalization of the KL property and KL functions is established. The sufficient conditions of the KL property are formally proved. The uniformized KL property is derived through extensions on compact sets, enabling convergence analysis in nonsmooth and nonconvex optimization problems.
    \item \textbf{Algorithm frameworks and convergence analysis.} i) A formal framework is provided for the proximal alternating linearized minimization in two-block problems \cite{bolte2014proximal}. Formalized convergence theorems such as critical point analysis through subdifferential characterizations and finite-length convergence under the KL property are provided. ii) A type of the ADMM algorithm is formalized for composite convex optimization problems with a linear constraint connecting two distinct components according to \cite{Fazel2013hankelmatrix, chen2017note}. Essential assumptions, including lower semicontinuity, convexity, and solvability of subproblems, are formally established to ensure systematic convergence analysis. Moreover, the error reduction bounds are formally presented, and global convergence of the ADMM method to KKT points is proved.
\end{enumerate}

The rest of this paper is organized as follows. The formalization of the subdifferential of general functions is given in Section \ref{sec: subdifferential}. KL property and KL functions are formalized in Section \ref{sec: KL}. We discuss the formalization of the scheme and the convergence analysis of the BCD and ADMM algorithms in Section \ref{sec: BCD} and Section \ref{sec: ADMM}, respectively. 

\section{Subdifferential Structure of General Function}
\label{sec: subdifferential}
\subsection{Definitions}
In this subsection, we review the definition and properties of subdifferential and the first-order optimality condition for general functions. These are particularly useful for analyzing nondifferentiable and nonconvex functions. Recall that for a convex function $f$ defined on a Hilbert space $E$, the subgradient $u$ provides a global lower bound for $f$. For all $y \in E$, the inequality $ f(y)\geq f(x)+\langle u,y-x\rangle$ holds at every point in the domain. However, when $f$ is nonconvex, this property is only expected to hold in a limiting sense. The definition of subdifferential is as follows.
\begin{definition}\label{def: F subdiff}
For a function $f : E  \to \mathbb{R}$ defined on a Hilbert space $E$, and a given point $x \in E$, the \emph{Fréchet subdifferential} is defined as
\begin{align*}
    \hat{\partial} f(x) = \left\{u\in\mathbb{R}^n \mid \liminf_{y \to x, y \neq x} \frac{f(y) - f(x) - \langle u, y - x \rangle}{\|y - x\|} \geq 0\right\}.
\end{align*}
\end{definition}
For convenience, we introduce the notion of \textit{f-attentive convergence}. The term $x^n \underset{f}{\to} x$ denotes that a sequence $\{x^n\}$ satisfies that $x^n \to x $ with $f(x^n) \to f(x)$ as $n \to \infty$. The limiting subdifferential is defined as follows.
\begin{definition}\label{def: limit subdiff}
The \emph{limiting subdifferential} of $f$ at $x\in E$ is defined as
    \begin{equation}
    \label{def: limit subdiff1}
\partial f(x) = \left\{u \in \mathbb{R}^n  \mid\exists x^k\underset{f}{\to} x, u^k \in \hat{\partial} f(x^k)\to u, \text{ as } k \to \infty\right\},
    \end{equation}
where $u$ is called a (general) subgradient of $f$ at $x$.
\end{definition}

Additional assumptions on the domain space are required to formalize these definitions. We assume the following instances for the space variable from this subsection.
\begin{lstlisting}
variable {E : Type*} [NormedAddCommGroup E] [InnerProductSpace ℝ E]
\end{lstlisting}
The assumptions in brackets give the properties of the space $E$, claiming that \lean{E} is a normed additive commutative group and an inner product space on $\mathbb{R}$. These assumptions can be automatically checked when using the related theorems. For example, when discussing functions on $\mathbb{R}^n$, Lean4 identifies that $\mathbb{R}^n$ is a normed add commutative group and an inner product space on $\mathbb{R}$. Hence, the theorems proved for general space $E$ work well for concrete space such as $\mathbb{R}^n$. 

The Fréchet subdifferential and limiting subdifferential in Lean are defined as follows.
\begin{lstlisting}
def differential_fun (x : E) (f : E → ℝ) (u : E) :=
    fun y ↦ Real.toEReal ((f y - f x - inner u (y - x)) / ‖ y - x ‖)

def f_subdifferential (f : E → ℝ) (x : E) : Set E :=
    {v : E | liminf (differential_fun x f v) (nhds[≠] x) ≥ 0}

def subdifferential (f : E → ℝ) (x : E) : Set E := 
    {v₀ : E | ∃ u : ℕ → E, Tendsto u atTop (nhds x) ∧ Tendsto (fun n ↦ f (u n)) atTop (nhds (f x)) ∧ (∃ v : ℕ → E, ∀ n, v n ∈ f_subdifferential f (u n) ∧ Tendsto v atTop (nhds v₀))}
\end{lstlisting}
The function \lean{Real.toEReal} is introduced to convert real numbers into extended real numbers (\lean{EReal}). This conversion is necessary for correctly applying the operator $\liminf$ in the definition of \lean{f_subdifferential}. Specifically, the limit infimum can attain infinite values only within extended real numbers. In contrast, within Lean's real number type ($\mathbb{R}$), the limit infimum defaults to zero when it would otherwise take an infinite value.

We can define the active domain of a function $f(x)$ using subdifferential, which consists of points whose subdifferential is not empty. The formalized definition is given as below.
\begin{lstlisting}
def active_domain (f : E → ℝ) : Set E := {x | subdifferential f x ≠ ∅}
\end{lstlisting}

In the above definitions, the regular subgradient uses the concept of \lean{liminf} in Lean. We can get an equal description using little-o notation through mathematical reformulation:
\begin{align*}
    f(y) \geq f(x) + \langle u, y - x \rangle + o(\|y - x\|), \quad \text{as}\; \|y-x\| \to 0.
\end{align*}
We can reformulate the definition using the concept of a filter to avoid explicitly managing Taylor expansions or infinitesimals, which brings convenience in some cases. An equivalent expression for the Fréchet subdifferential is defined in the following statement.
\begin{theorem}\label{thm: equvilant subdifferential}
For a given $x \in \operatorname{dom} f$, $u \in \hat{\partial} f(x) $ if and only if for any $ \varepsilon > 0$ and every $y$ in the neighborhood of x, it holds that
\begin{align*}
    f(y) - f(x) - \langle u, y - x \rangle \geq -\varepsilon \|y - x\|.
\end{align*}
\end{theorem}
The above theorem is formalized as follows.
\begin{lstlisting}
theorem has_f_subdiff_iff : u ∈ f_subdifferential f x ↔ ∀ ε > 0, 
(nhds x).Eventually fun x => f y - f x - inner u (y - x) ≥ -ε * ‖ y - x ‖  
\end{lstlisting}
To state the next theorem, we define:
\begin{align}
G (f,\hat{\partial} f) & := \left\{ (x,f(x),u)|x \in \mathbb{R}^n,u \in \hat{\partial} f(x)\right\},\\
G (\partial f) & := \left\{ (x,u)|x \in \mathbb{R}^n,u \in \partial f(x)\right\}.
\end{align}
We next introduce a more streamlined description of the subdifferential by leveraging the closure of the subdifferential graph. This reformulation provides an intuitive geometric interpretation by directly relating the subdifferential to the closure of its graph. The following theorem states the closedness of the subdifferential graph. This property is also used in the proof of the KL property.
\begin{theorem}
    Let $(x_k,x_k^*)\in G (\partial f)$ be a sequence that converges to $(x,x^*)$.
    If $f(x_k)$converges to $f(x)$, then $(x,x^*)\in G(\partial f)$. Or equivalently,
    \begin{align}
        \partial f(x) = \left\{ u \in \mathbb{R}^n | (x, f(x), u) \in \overline{G (f, \hat{\partial} f)}  \right\}. \label{eq: subdifferential definition 2}
    \end{align}
\end{theorem}
The theorem is formalized as below.
\begin{lstlisting}
def subdifferential_Graph (f : E → ℝ) := 
{(x, u) : E × E | u ∈ subdifferential f x}

theorem GraphOfSubgradientIsClosed {f : E → ℝ} {xn un : ℕ → E} {x u : E}
(hx : ∀ n , (xn n, un n) ∈ subdifferential_Graph f)
(hconv : Tendsto (fun n => (xn n , un n)) atTop (nhds (x, u)))
(hf : Tendsto (fun n => f (xn n)) atTop (nhds (f x))) :
(x, u) ∈ subdifferential_Graph f 
\end{lstlisting}

\subsection{Properties of Subdifferentials}
These definitions of subdifferentials enable the derivation of the first-order optimality condition of general functions. For general functions, the first-order condition is stated as follows.
\begin{theorem}
If $x^*$ is a local minimum point of $f$, then $0 \in \partial f (x^*)$.
\label{thm: optimal}
\end{theorem}
In Lean, we formalize this as follows.
\begin{lstlisting}
theorem first_order_optimality_condition (f : E → ℝ) (x₀ : E) 
(hx : IsLocalMin f x₀) : (0 : E) ∈ f_subdifferential f x₀
\end{lstlisting}
The point \( y \) satisfying \( 0 \in \partial f(y) \) is called a critical point of the function \( f(x) \). The set of critical points of $f(x)$ is denoted as $\operatorname{crit }f (x)$. This is formalized as follows.
\begin{lstlisting}
def critial_point (f : E → ℝ) : Set E := {x | 0 ∈ subdifferential f x}
\end{lstlisting}
It is evident that a local minimum must be a critical point of the objective function. 

Notably, both the Fréchet subdifferential and the limiting subdifferential reduce to the gradient of $f$ at points where the function $f$ is smooth. Several theorems will be introduced to state the equivalence between subdifferential and gradient for differentiable functions. Since the gradient of the function is considered, we additionally assume \lean{[CompleteSpace E]} to align with the assumptions regarding the existence of the gradient. To begin with, we present the following theorem to show the definition of Fréchet differential is the same as the gradient when the function is differentiable.
\begin{theorem}\label{thm: subdifferential diff}
    The following two statements are equivalent:
    \begin{enumerate}
        \item $f$ is differentiable at $x$ with a gradient $u$. 
        \item $u \in \hat{\partial}f(x)$ and $-u \in \hat{\partial}(-f)(x)$.
    \end{enumerate}
\end{theorem}
The formalized version can be given as below.
\begin{lstlisting}
theorem HasGradientAt_iff_f_subdiff : HasGradientAt f u x ↔ u ∈ f_subdifferential f x ∧ -u ∈ f_subdifferential (-f) x
\end{lstlisting}
The proof of this proposition mainly utilizes the equivalent transformation of Fréchet subdifferential provided by \lean{has_f_subdiff_iff} in Theorem \ref{thm: equvilant subdifferential}, thereby avoiding dealing with the lower limit in the definition of Fréchet subdifferential. The following property states the equivalence between gradient, Fréchet subdifferential, and limiting subdifferential.
\begin{theorem}\label{thm: singleton subdifferential}
    Let $f(x)$ be a differentiable function with gradient $\nabla f(x)$ on a Hilbert space $E$. Then for every $z \in E$ , $\hat{\partial}f(z)$ is a singleton set with the element $\nabla f(z)$ . 
\end{theorem}
The formalization of the theorem is given as follows.
\begin{lstlisting}
theorem f_subdiff_gradient (f' : E → E) (z : E) (hf : HasGradientAt f (f' z) z) : f_subdifferential f z = {f' z}
\end{lstlisting}
We prove this base on Theorem \ref{thm: subdifferential diff}. The conclusion $\nabla f(z) \in \hat{\partial}f(z)$ is a direct result from the proposition \lean{HasGradientAt_iff_f_subdiff}. The theorem \lean{f_subdiff_neg_f_subdiff_unique} is introduced to establish the connection between the subdifferentials of functions that are negatively correlated. Theorem \ref{thm: singleton subdifferential} can be deduced from these preparations. 

The calculation rules of the subdifferential are also formalized. The following theorems show that under certain conditions, the Fréchet subdifferential of the sum of two functions can be converted into the sum of the Fréchet subdifferentials of the two functions separately.
\begin{theorem}\label{thm: f_subdifferential rule}
    Let $f(x), g(x)$ be two functions defined on $E$. If $g(x)$ is differentiable, then it holds that $\hat{\partial}(f + g)(x) = \hat{\partial}f(x) + \hat{\partial}g(x)$. 
\end{theorem}
The formalized version of the theorem is as follows.
\begin{lstlisting}
theorem f_subdiff_add (f g : E → ℝ) (g' : E → E) (x : E) (hg : ∀ (x₁ : E), HasGradientAt g (g' x₁) x₁) (z : E) :
    z ∈ f_subdifferential (f + g) x ↔ ∃ zf ∈ f_subdifferential f x, ∃ zg ∈ f_subdifferential g x, z = zf + zg
\end{lstlisting}
With an extra condition that $g$ is continuously differentiable, Theorem \ref{thm: f_subdifferential rule} also holds for the limiting subdifferential.
\begin{theorem}
    Let $f, g$ be two functions defined on $E$. If $g$ is continuous differentiable on $E$, it holds that ${\partial}(f + g) = {\partial}f + {\partial}g$.
\end{theorem}
The formalized theorem gives as follows.
\begin{lstlisting}
theorem subdiff_add (f g : E → ℝ) (g' : E → E) (x : E) (hg : ∀ (x₁ : E), HasGradientAt g (g' x₁) x₁) (gradcon : Continuous g') (z : E) : z ∈ subdifferential (f + g) x ↔ ∃ zf ∈ subdifferential f x, ∃ zg ∈ subdifferential g x, z = zf + zg
\end{lstlisting}

The proximal operator is also widely used in optimization algorithms, which addresses an optimization problem that balances maintaining proximity to a specified input point $x$ while also minimizing the function $f$. For general functions, the proximal operator maps a point to a set as follows.
\begin{definition}
The \emph{proximal operator} for function $f$ at point $x$ is given as
\begin{align}\label{eq: prox_def}
    \text{prox}_f(x) = \underset{u}{\arg\min} \left\{ f(u) + \frac{1}{2}\|u-x\|^2 \right\}.
\end{align}
\end{definition}
The formalization of the proximal operator can be referred to \cite{li2024formalization}. The optimality condition of the optimization problem naturally links the proximal operator with the subdifferential. The connection is described in the following theorem. 
\begin{theorem}
If $f(x)$ is a function from a Hilbert space $E$ to $\mathbb{R}$, it holds that
\begin{align}
    u \in \text{prox}_f(x) \quad \text{if and only if} \quad x - u \in \partial f(u).
    \label{eq: proximal connection}
\end{align}
\end{theorem}
The theorem is formalized as below.
\begin{lstlisting}
theorem f_subdiff_prox {f : E → ℝ} {x y : E} (h : prox_prop f y x) : y - x ∈ f_subdifferential f x 
\end{lstlisting}

A key step to prove the above propositions is to utilize first order optimality condition proven before to connect the minimality definition of the proximal operator and Fréchet subdifferential. The equivalence and operational properties of Fréchet subdifferential for certain functions also help to simplify the problem.

\section{The KL Property and KL Function}
\label{sec: KL}
In the past few decades, the KL property has emerged as a versatile analytical framework for investigating convergence rates of first-order optimization methods, particularly when dealing with objectives that exhibit both nonsmoothness and nonconvexity, see \cite{attouch2013convergence, bolte2014proximal, bento2024convergence}. It is an extension to the Polyak-Łojasiewicz (PL) inequality developed by Polyak \cite{polyak1963gradient} and Łojasiewicz \cite{lojasiewicz1963propriete},
\begin{equation*}
    \|\nabla f(x) \| \geq \mu \| f(x) - f(\bar x)\|^{\frac{1}{2}}, \quad \mu > 0,
\end{equation*}
where $f$ is an $L$-smooth function, and $\bar x$ is a global minimum of $f$. The linear convergence of the gradient descent algorithm can be established by this well-known inequality. It illustrates that the sub-optimality can be controlled by the vanishing of gradient norm and a special concave function $\phi(t) = \mu t^{\frac{1}{2}}$. Such functions are the so-called desingularizing functions in the context of the KL property.
\begin{definition}
\label{def : desingularizing}
    For any $\eta \in (0,+\infty)$, if $\phi$ satisfies the following properties:
    \begin{enumerate}
        \renewcommand{\theenumi}{\roman{enumi}} 
        \item $\phi:[0,\eta) \to [0, \infty)$ is continuous and concave;
        \item $\phi(0) = 0$;
        \item $\phi$ is continuously differentiable on $(0,\eta)$ with $\phi'>0$ on $(0,\eta)$,
    \end{enumerate}
    then we call $\phi$ a desingularizing function and denote it as $\phi \in \Phi_\eta$.
\end{definition}
The corresponding formalization of this family of functions is as follows.
\begin{lstlisting}
def desingularizing_function (η : ℝ) := {φ : ℝ → ℝ | (ConcaveOn ℝ (Ico 0 η) φ) ∧ (φ 0 = 0) ∧ (ContDiffOn ℝ 1 φ (Ioo 0 η)) ∧ (ContinuousAt φ 0) ∧ (∀ x ∈ Ioo 0 η, deriv φ x > 0)}
\end{lstlisting}

These desingularizing functions are essential to quantify the behavior of the objective $f$ near its critical points, enabling the analysis of the relationship between the subdifferential and the function values. We can define the KL property utilizing such functions. Still, all variables are considered in the Hilbert space $E$.
\begin{definition}
\label{def:KL}
    Let $f: E \to \mathbb{R}$ be a lower semicontinuous function. We say that $f$ satisfies the KL property at $u\in \mathrm{dom}\partial f := \{ u \in E : \partial f (u) \neq \emptyset \}$ if there exist $\eta \in (0, \infty)$, a neighborhood $S$ of $u$ and $\phi\in \Phi_\eta$ such that for any $x \in S \cap [f(u) < f < f(u) + \eta]$, the following inequality holds:
        \begin{equation*}
            \phi'(f(x)-f(u))\cdot \mathrm{dist}(0, \partial f(x)) \geq 1.
        \end{equation*}
    The point $u$ is called a KL point of $f$. If $f$ satisfies KL property at every point of $\mathrm{dom} \partial f$, we say that $f$ is a KL function.
\end{definition}
We generalize the assumption of lower semicontinuity in Definition \ref{def:KL} in formalization. Note that the KL property holds naturally if $\partial f (\bar u) = \emptyset$. Hence, the non-emptiness of subdifferential can also be omitted in formalization. These two assumptions need not be included in the formalization of the KL point and KL function. Consequently, $\mathrm{dist}(0,\partial f(u))$ might be infinity when $\partial f(u)$ is empty. We use \lean{EMetric.infEdist} to keep the inequality on extended non-negative real numbers \lean{ENNReal}, and \lean{ENNReal.ofReal} to convert real numbers to \lean{ENNReal}. The formalization of the KL point and KL function are as follows.
\begin{lstlisting}
def KL_point (f : E → ℝ) (u : E) : Prop := ∃ η ∈ Ioi 0, ∃ s ∈ nhds u, ∃ φ ∈ desingularizing_function η, ∀ x ∈ s ∩ {y | f u < f y ∧ f y < f u + η}, (ENNReal.ofReal (deriv φ (f x - f u))) * (EMetric.infEdist 0 (subdifferential f x)) ≥ ENNReal.ofReal 1

def KL_function (f : E → ℝ) : Prop := ∀ u ∈ active_domain f, KL_point f u
\end{lstlisting}

It is shown in \cite{attouch2010proximal} that the KL property holds at noncritical points of $f$.
\begin{theorem}
    Suppose that $f:E \to \mathbb{R}$ is a lower semicontinuous function. The KL property holds at $x$, if it holds that $0 \notin \partial f(x)$.
\end{theorem}
This theorem is formalized as follows.
\begin{lstlisting}
theorem KL_property_at_noncritical_point (h_noncrit : 0 ∉ subdifferential f x) : KL_point f x 
\end{lstlisting}

The uniformized KL property \cite{bolte2014proximal} is incorporated to deal with the case when the KL property holds at many different points.

\begin{theorem}
    Suppose $\Omega$ is a compact set and  $f:E \to \mathbb{R}$ is a lower semicontinuous function. Suppose $f$ is constant on $\Omega$ and satisfies the KL property at every point of $\Omega$. Then there exist $\epsilon>0, \eta>0$, and a function $\phi \in \Phi_\eta$, such that for all $u\in \Omega$ and all $x$ satisfying $\mathrm{dist}(x,\Omega) < \epsilon, f(u) < f (x) < f (u)+ \eta$, it holds
    \begin{equation*}
        \phi'(f(x) - f(u)) \cdot \mathrm{dist}(0, \partial f (x)) \geq 1.
    \end{equation*}
\end{theorem}
\begin{lstlisting}
theorem uniformized_KL_property {f : E → ℝ} {Ω : Set E} (h_compact : IsCompact Ω) (h_Ω1 : ∀ x ∈ Ω, KL_point f x) (h_Ω2 : is_constant_on f Ω) : ∃ ε ∈ Ioi 0, ∃ η ∈ Ioi 0, ∃ φ ∈ desingularizing_function η, ∀ u ∈ Ω , ∀ x ∈ {y : E | (EMetric.infEdist y Ω).toReal < ε} ∩ {y | f u < f y ∧ f y < f u + η}, (ENNReal.ofReal (deriv φ (f x - f u))) * EMetric.infEdist 0 (subdifferential f x) ≥ 1
\end{lstlisting}
To prove this theorem, a finite open cover can be chosen based on the neighborhoods of the KL points using \lean{IsCompact.elim_finite_subcover_image}. Then there exists an open ball in this finite cover from \lean{IsCompact.exists_thickening_subset_open}. Furthermore, one can prove that the sum of the finite desingularizing functions on the finite cover will be in $\Phi_\eta$ for some $\eta > 0$ by checking Definition \ref{def : desingularizing}. Finally, a uniformized neighborhood and a desingularizing function are constructed for the inequality.

\section{The Block Coordinate Descent Algorithm}
\label{sec: BCD}
In this section, we consider a special form of the BCD algorithm with proximal alternating linearized minimization following \cite{bolte2014proximal}. The BCD algorithm is applied to two-block optimization problems as:
\begin{align}
    \underset{x, y}{\min} \; \Psi(x,y) = f(x) + g(y) + H(x,y),
    \label{eq : target BCD}
\end{align}
where $f(x) : \mathbb{R}^n \to \mathbb{R}$ and $g(y) : \mathbb{R}^m \to \mathbb{R}$ are separated parts for $x$ and $y$, respectively. $H(x,y) : \mathbb{R}^n \times \mathbb{R}^m \to \mathbb{R}$ is the composite part. The basic assumptions for this problem are listed as below.
\begin{assumption}\label{assumption: BCD}
 \begin{enumerate}
        \renewcommand{\theenumi}{\roman{enumi}} 
        \item The separated parts $f(x)$ and $g(y)$ are lower semicontinuous functions. $f(x)$, $g(y)$ and $\Psi(x,y)$ are all bounded below.
        \item $H(x,y)$ is a differentiable function with a Lipschitz continuous gradient. The Lipschitz constant of the gradient is $l$, i.e.
        \begin{align*}
            \left\| \left( \nabla_x H(x_1, y_1) - \nabla_x H(x_2, y_2), \nabla_y H(x_1, y_1) - \nabla_y H(x_2, y_2) \right) \right\| \leq l \left\| (x_1 - x_2, y_1 - y_2) \right\|.
        \end{align*}
    \end{enumerate}
\end{assumption}
The BCD algorithm utilizes an alternating optimization method to update $x$ and $y$, where each update of $y$ integrates the latest value of $x$. In each iteration, the linear approximation of the smooth part $H(x,y)$ is used to approximate the local property of $H(x,y)$. Another penalty term with step size is added to balance the error of approximation. For the $k$-th iteration, we solve the following subproblems:
\begin{align}
\begin{aligned}
    x^{k+1} &\in \underset{x}{\arg\min} \left\{f(x)+ \langle \nabla_x H(x^k, y^k), x - x ^k \rangle + \frac{c_k}{2} \|x-x^k\|^2  \right\}, \\
    y^{k+1} &\in \underset{y}{\arg\min}\left\{g(y) + \langle \nabla_y H(x^{k+1}, y^k), y - y^k \rangle + \frac{d_k}{2} \|y-y^k\|^2  \right\}.
\end{aligned}
\label{eq : BCD update ori}
\end{align}
The set of minimizers of the subproblem may not be a singleton. We select an arbitrary element to continue the iterations. The parameters $c_k$ and $d_k$ control the step sizes. Through the definition of the proximal operator \eqref{eq: prox_def}, the subproblems can be rewritten as
\begin{align}
\begin{aligned}
    x^{k+1} &\in \text{prox}_{c_k f} \left( x^k - c_k \nabla_x H(x^k, y^k) \right), \\
    y^{k+1} &\in \text{prox}_{d_k g} \left( y^k - d_k \nabla_y H(x^{k+1}, y^k) \right).
\end{aligned}
\label{eq : BCD update}
\end{align}

\subsection{Formalization of the Update Scheme}
For any numerical algorithm, the formalization can typically be decomposed into two primary components: the problem data and the update mechanism. The problem data encapsulate well-founded assumptions about the problem, whereas the update mechanism encompasses the initial point, sequence of updates, and the specific rules applied at each iteration. The formalized problem data is presented below:
\begin{lstlisting}
variable [NormedAddCommGroup E] [InnerProductSpace ℝ E] [CompleteSpace E]
variable [NormedAddCommGroup F] [InnerProductSpace ℝ F] [CompleteSpace F]
structure ProblemData (f : E → ℝ) (g : F → ℝ) (H : (WithLp 2 (E × F)) → ℝ) (l : NNReal) : Prop where
  lbdf : BddBelow (f '' univ)     lbdg : BddBelow (g '' univ)
  hf : LowerSemicontinuous f     hg : LowerSemicontinuous g
  conf : ContDiff ℝ 1 H          lpos : l > (0 : ℝ)
  lip : LipschitzWith l (gradient H)
\end{lstlisting}
All assumptions are bundled into the structure problem data. This facilitates the identification and use of the problem conditions. The concrete space $\mathbb{R}^n$ is replaced by a Hilbert space \lean{E}. This not only makes the analysis clearer in formalization, but also generalizes the current results on $\mathbb{R}^n$ to a general Hilbert space.  It is essential to note that the function \(H\) operates in the domain \lean{WithLp 2 (E × F)}. Handling the product space \((x, y)\) is challenging. Within \mathlib, two types of product spaces are defined: \lean{E × F} and \lean{WithLp 2 (E × F)}. Although the elements of these spaces can both be represented as \((x, y)\), their norms differ significantly. For \lean{(x, y)} of type \lean{E × F}, the norm is defined as \(\|(x, y)\| = \max \{\|x\|, \|y\|\}\). However, for \lean{WithLp 2 (E × F)}, it holds that \(\|(x, y)\| = \sqrt{\|x\|^2 + \|y\|^2}\). The former space is merely normed and lacks an inner product structure, while the latter represents the typical inner product space used for product spaces. These norms are equivalent, each bounded by absolute constants. Since the usage of gradient requires that the domain is a Hilbert space, \lean{WithLp 2 (E × F)} is a better choice choice for defining gradients due to its inner product space properties.

We give the formalized algorithm towards \eqref{eq : BCD update} using the \lean{structure} type extending over \lean{ProblemData}. The update scheme, predicated on the problem data, is structured as follows:
\begin{lstlisting}
structure BCD (f : E → ℝ) (g : F → ℝ) (H : (WithLp 2 (E × F)) → ℝ) (l : NNReal)(x0 : E) (y0 : F) extends ProblemData f g H l where
  (x : ℕ → E)  (y : ℕ → F)
  (x0 : x 0 = x0) (y0 : y 0 = y0) (c d : ℕ → ℝ)
  s₁ : ∀ k, prox_prop (c k • f) (x k - c k • (grad_fst H (y k) (x k))) (x (k + 1))
  s₂ : ∀ k, prox_prop (d k • g) (y k - d k • (grad_snd H (x (k + 1)) (y k))) (y (k + 1))
\end{lstlisting}
The definitions for \(f\), \(g\), and \(H\) are directly inherited from the problem data. The term \lean{prox_prop f x y} means that $y \in \operatorname{prox}_f(x)$. We introduce additional notations using the same namespace BCD as follows:
\begin{lstlisting}
def BCD.z {self : BCD f g H l x0 y0} : ℕ → WithLp 2 (E × F) :=
    fun n ↦ (WithLp.equiv 2 (E × F)).symm (self.x n, self.y n)

def BCD.ψ {_ : BCD f g H l x0 y0} := 
    fun z : WithLp 2 (E × F) ↦ (f z.1 + g z.2 + H z)
\end{lstlisting}
These definitions are given in the local namespace as beginning with ``\lean{BCD.}”. Hence, if we use \lean{alg : BCD f g H l x0 y0} to define an algorithm, \lean{alg.z} automatically turns to this definition, which gives convenience for complicated calculation.

\subsection{Convergence Analysis}

For the analysis of the fixed step size version of BCD, we mainly follow \cite{bolte2014proximal} and the detailed proof is omitted in this paper. We set the step sizes \(c_k\) and \(d_k\) consistently as \(\frac{1}{\gamma l}\), where \(\gamma\) serves as a critical hyper-parameter. The parameter \(l\) represents the Lipschitz constant of \(\nabla H\). The conditions on the stepsizes are given in \lean{variable} as follows and for simplicity we omit to claim them for each theorem in this subsection.
\begin{lstlisting}
variable (γ : ℝ) (hγ : γ > 1) (ck: ∀ k, alg.c k = 1 / (γ * l)) (dk: ∀ k, alg.d k = 1 / (γ * l))
\end{lstlisting}

To demonstrate the convergence of the algorithm, the analysis unfolds in several stages. Initially, we establish a sufficient descent condition showing that the function values decremented by the algorithm are monotonically decreasing.
\begin{theorem}
\label{thm: sufficient descent bcd}
Suppose that Assumption \ref{assumption: BCD} holds. Under the condition that the iterative sequence $\{z^k = (x^k, y^k)\}$ generated by updates \eqref{eq : BCD update} is bounded and the step sizes satisfy $c_k = d_k = \frac{1}{\gamma l}$ for $\gamma > 1$, the sequence exhibits the following characteristics.
\begin{enumerate}
    \item At each iteration, the decrement property holds:
   \[
   (\gamma - 1)l \|z^{k+1} - z^k\|^2 \leq H(z^k) - H(z^{k+1}), \quad \forall k \geq 0.
   \]
   \item The sequence $\{\|z^{k+1} - z^k\|\}_{k=1}^\infty$ is square summable:
    \[
    \sum_{k=1}^\infty \|z^{k+1} - z^k\|^2 < \infty,
    \]
    ensuring that $\lim\limits_{k \to \infty} \|z^{k+1} - z^k\| = 0.$
\end{enumerate}
\end{theorem}
The following two statements formalize Theorem \ref{thm: sufficient descent bcd}.
\begin{lstlisting}
theorem Sufficient_Descent1 {alg : BCD f g H l x0 y0} :
    ∃ ρ₁ > 0, ρ₁ = (γ - 1) * l ∧ ∀ k, ρ₁ / 2 * ‖ alg.z (k+1) - alg.z k ‖ ^ 2 ≤ alg.ψ (alg.z k) - alg.ψ (alg.z (k + 1)) 
\end{lstlisting}
\begin{lstlisting}
theorem Sufficient_Descent2 {alg : BCD f g H l x0 y0} (lbdψ : BddBelow (alg.ψ '' univ)):
    Tendsto (fun k ↦ ‖ alg.z (k + 1) - alg.z k ‖) atTop (nhds 0) 
\end{lstlisting}
Theorem \lean{Sufficient_Descent1} establishes the non-increasing nature of the objective function. Theorem \lean{Sufficient_Descent2} guarantees that the distances between successive updates converge to zero, ensuring update stabilization. 

To gain a deeper understanding of the analysis, it is crucial to examine the behavior of the subdifferential at the update points. We can establish that the subdifferential at each iteration is upper-bounded by the norm of the difference between two consecutive points, which implies that the norm of the subdifferential tends to zero as the iterations progress.
\begin{lemma}\label{lemma: BCD subdifferential}
Suppose that Assumptions \ref{assumption: BCD} holds. Suppose $\{z^k\}$ generated by scheme \eqref{eq : BCD update} is bounded and the step sizes satisfy $c_k = d_k = \frac{1}{\gamma l}$ for $\gamma > 1$. For each positive integer k, define
\begin{align*}
    A_x^k := c_{k-1} (x^{k-1} - x^k) + \nabla_x H(x^k, y^k) - \nabla_x H(x^{k-1}, y^k), \\
    A_y^k := d_{k-1} (y^{k-1} - y^k) + \nabla_y H(x^k, y^k) - \nabla_y H(x^k, y^{k-1}).
\end{align*}
Then  \( (A_x^k, A_y^k) \in \hat{\partial} \Psi(x^k, y^k) \) and there exists \( M = (2\gamma+2)l > 0 \) such that
\begin{align*}
    \left\| (A_x^k, A_y^k) \right\| \leq \left\| A_x^k \right\| + \left\| A_y^k \right\| \leq M \left\| z^k - z^{k-1} \right\|, \quad \forall k \geq 1, 
\end{align*}
\end{lemma}
The specific values of $A_x^k$ and $A_y^k$ for the Fréchet subdifferentials are not essential to the continuation of the proof. What matters is their existence and norm bound. Therefore, for simplicity, our formalized version of Lemma \ref{lemma: BCD subdifferential} omits the explicit values of the subdifferentials and focuses solely on their existence and bound. In our proof, the values of $A_x^k$ and $A_y^k$ are used, and the membership relation is established by the calculation rule in Theorem \ref{thm: f_subdifferential rule}. One of the most difficult problem in formalization is to deal with the subdifferential with the product space. We establish calculation rules for block-wise subdifferential. The formalized version is provided below.
\begin{lstlisting}
theorem Ψ_subdiff_bound {alg : BCD f g H l x0 y0} :
    ∃ ρ > 0, ∀ k, ∃ dΨ ∈ f_subdifferential alg.ψ (alg.z (k + 1)), ‖dΨ‖ ≤ ρ * ‖ alg.z (k + 1) - alg.z k ‖ 
\end{lstlisting}
From this boundedness of the subdifferential, we obtain that the limit points of the update sequence \(\{z^k\}\) constitute critical points of the objective function.
\begin{lemma}\label{lemma: BCD critical}
Suppose that Assumptions \ref{assumption: BCD} holds. Let $\{z^k\}$ generated by algorithm \eqref{eq : BCD update} be bounded and the step sizes satisfy $c_k = d_k = \frac{1}{\gamma l}$ for $\gamma > 1$. Denote the limit set of $\{z^k\}$ as $\omega(z^0)$. It holds that $\emptyset \neq \omega(z^0) \subset \operatorname{crit} \Psi(z)$ and the objective function $\Psi(z)$ is finite and constant on $\omega(z^0)$.
\end{lemma}
We add a definition of limit set in Lean using the function \lean{MapClusterPt} in Mathlib4, which denotes the cluster point of a mapping. The formalized theorem is given as below.
\begin{lstlisting}
def limit_set (z : ℕ → E) := {x | MapClusterPt x atTop z}
theorem limitset_property {alg : BCD f g H l x0 y0} (bd : Bornology.IsBounded (alg.z '' univ)) (lbdψ : BddBelow (alg.ψ '' univ)):
    (limit_set alg.z).Nonempty ∧ (limit_set alg.z) ⊆ critial_point alg.ψ ∧ ∃ (c : ℝ) , ∀ x ∈ (limit_set alg.z) , alg.ψ x = c 
\end{lstlisting}

Theorem \ref{thm: sufficient descent bcd} does not address the properties of the limit point of \( z^k \). To ensure convergence to a critical point of the objective function \( \Psi(x, y) \), it is assumed that the objective function \(\Psi(x,y)\) exhibits the KL property near these critical points, strengthening conditions around the critical points of \(H(x, y)\). We can obtain a stronger version of the convergence theorem. 
\begin{theorem}
\label{thm: bounded length bcd}
Suppose that Assumption \ref{assumption: BCD} holds. If $\Psi(x,y)$ is a KL function and the iterative sequence $\{z^k\}$ adheres to the updates \eqref{eq : BCD update}, the following theorems hold.
\begin{enumerate}
    \item The sequence $\{z^k\}$ possesses a finite total variation:
   \[
   \sum_{k=1}^\infty \|z^{k+1} - z^k\| < +\infty,
   \]
   \item $\{z^k\}$ converges to a critical point $z^* = (x^*, y^*)$ of $\Psi$.
\end{enumerate}
\end{theorem}
We can prove the formalized version of Theorem $\ref{thm: bounded length bcd}$.
\begin{lstlisting}
theorem Limited_length {alg : BCD f g H l x0 y0} (bd : Bornology.IsBounded (alg.z '' univ)) (hψ : KL_function alg.ψ) (lbdψ : BddBelow (alg.ψ '' univ)): ∃ M : ℝ, ∀ n,
    ∑ k in Finset.range n, ‖alg.z (k + 1) - alg.z k‖ ≤ M 
\end{lstlisting}
\begin{lstlisting}
theorem Convergence_to_critpt {alg : BCD f g H l x0 y0} (bd : Bornology.IsBounded (alg.z '' univ)) (hψ : KL_function alg.ψ) (lbdψ : BddBelow (alg.ψ '' univ)) : 
    ∃ z_ : (WithLp 2 (E × F)), z_ ∈ (critial_point alg.ψ) ∧ Tendsto alg.z atTop (nhds z_)
\end{lstlisting}
The proofs of the aforementioned theorems rely on algebraic manipulations and basic analytical tools, including inequalities and continuity arguments, to establish the sufficiency and boundedness conditions. This rigorous framework supports a clear iterative convergence analysis and enables precise customization of the BCD algorithm for specific problem settings. Furthermore, the foundational structure naturally extends to a range of optimization challenges.

\section{The Alternating Direction Method of Multipliers}
\label{sec: ADMM}
The ADMM method is frequently used to solve optimization problems comprising multiple components coupled by a linear constraint. In this paper, we focus on a certain kind of the ADMM method:
\begin{equation}
    \begin{aligned}
        \min_{x_1, x_2} & \quad f_1(x_1) + f_2(x_2), \\
        \text{s.t.} & \quad A_1 x_1 + A_2 x_2 = b,
    \end{aligned}
    \label{eq: admm form}
\end{equation}
where $f_1(x_1) : \mathbb{R}^n \to \mathbb{R}$ and $f_2(x_2) : \mathbb{R}^m \to \mathbb{R}$ are two separated parts of the objective function. The matrices \( A_1 \in \mathbb{R}^{p \times n} \) and \( A_2 \in \mathbb{R}^{p \times m} \) link the variables via linear constraints, with \( b \in \mathbb{R}^p \) representing a given vector. Our analysis of this algorithm mainly follows \cite{Fazel2013hankelmatrix, chen2017note}. The assumptions of this problem are given below.
\begin{assumption}\label{assumption: ADMM}\
 \begin{enumerate}
        \renewcommand{\theenumi}{\roman{enumi}} 
        \item The separated parts $f_1(x_1)$ and $f_2(x_2)$ are lower semicontinuous convex functions. 
        \item The solution set of problem \eqref{eq: admm form} is nonempty and satisfies Slater's condition, i.e. there exists a feasible solution in the interior of the domain:
        \begin{align*}
            \exists (\tilde{x}_1, \tilde{x}_2) \in \text{int} ( \text{dom} f_1 \times \text{dom} f_2 ) \quad \text{s.t.} \quad A_1 \tilde{x}_1 + A_2 \tilde{x}_2 = b.
        \end{align*}
    \end{enumerate}
\end{assumption}
The augmented Lagrangian function of problem \eqref{eq: admm form} is given as:
\begin{align}
    L_{\rho}(x_1, x_2, y) = f_1(x_1) + f_2(x_2) + \left<y, A_1x_1 + A_2x_2 - b\right> + \frac {\rho} 2  \left\lVert A_1x_1 + A_2x_2 - b\right \rVert^2,
    \label{eq: augmented lagrangian}
\end{align}
where \( \rho \) is the penalty parameter. To solve the constraint optimization problem \eqref{eq: admm form}, the ADMM method sequentially minimizes the augmented Lagrangian function with respect to each variable, then updates the Lagrange multipliers. Based on \eqref{eq: augmented lagrangian}, the iterative scheme of ADMM is as follows:
\begin{align}
        x_{1}^{k+1} &= \underset{x_1}{\arg \min}\; L_{\rho}(x_1,x_2^k,y^k), \label{eq: admm update 1}\\
        x_{2}^{k+1} &= \underset{x_2}{\arg \min}\; L_{\rho}(x_1^{k+1},x_2,y^k), \label{eq: admm update 2}\\
        y^{k+1} &= y^k + \tau \rho (A_1x_1^{k+1} + A_2x_2^{k+1} - b), \label{eq: admm update 3}
\end{align}
where \( \tau \) is a step size. The alternating minimization steps for \( x_1 \) and \( x_2 \) are followed by the iteration for the dual variable \( y \). We have an additional assumption on the scheme \eqref{eq: admm update 1}-\eqref{eq: admm update 3} as below. 
\begin{assumption}\label{assumption: ADMM 2}
    The subproblems \eqref{eq: admm update 1} and \eqref{eq: admm update 2} have a unique solution, respectively.
\end{assumption}
The first item in Assumption \ref{assumption: ADMM} guarantees the convexity of the subproblems. Assumption \ref{assumption: ADMM 2} ensures the uniqueness of the solutions to the subproblems, which is crucial for well-defined iterations. This assumption, though seemingly strong, is necessary to ensure the convergence of iteration points to a KKT point of problem \eqref{eq: admm form}. Under the Slater’s condition, the solution to problem \eqref{eq: admm form} can be expressed as a KKT pair, enabling us to replace the solution with its corresponding KKT pair and simplify the analysis.

\subsection{Formalization of the Update Scheme}
This subsection provides a detailed formalization of the ADMM method. First, we define the fundamental structure of the optimization problem, given by \eqref{eq: admm form}. This structure is formally represented using a \lean{class} type in Lean as follows.
\begin{lstlisting}
class Opt_problem [FiniteDimensional ℝ E₁] [FiniteDimensional ℝ E₂] [FiniteDimensional ℝ F] where
  f₁ : E₁ → ℝ            f₂ : E₂ → ℝ
  A₁ : E₁ → L[ℝ] F       A₂ : E₂ → L[ℝ] F        
  b  : F
  lscf₁ : LowerSemicontinuous f₁        lscf₂ : LowerSemicontinuous f₂
  cf₁ : ConvexOn ℝ univ f₁              cf₂ : ConvexOn ℝ univ f₂
  nonempty : ∃ x₁ x₂, (A₁ x₁) + (A₂ x₂) - b = 0 ∧ IsMinOn (fun (x₁, x₂) ↦ f₁ x₁ + f₂ x₂) univ (x₁, x₂)
\end{lstlisting}
The functions $f_1$ and $f_2$ are real-valued, lower semicontinuous (\lean{lscf₁}, \lean{lscf₂}), and convex (\lean{cf₁}, \lean{cf₂}). The domains of the functions are taken as the universe for simplicity. The property \lean{nonempty} guarantees the second item in Assumption \ref{assumption: ADMM}. Although $A_1$ and $A_2$ are originally defined as matrices, they are implemented as finite-dimensional linear maps in the code to enhance usability. 
The following definition captures the augmented Lagrangian function \eqref{eq: augmented lagrangian}.

\begin{lstlisting}
def Augmented_Lagrangian_Function (opt : Opt_problem E₁ E₂ F) (ρ : ℝ) : E₁ × E₂ × F → ℝ :=  fun (x₁, x₂, y) => (opt.f₁ x₁) + (opt.f₂ x₂) + inner y ((opt.A₁ x₁) + (opt.A₂ x₂) - opt.b) + ρ / 2 * ‖(opt.A₁ x₁) + (opt.A₂ x₂) - opt.b ‖ ^ 2
\end{lstlisting}

ADMM alternately updates $x_1$ and $x_2$, minimizing one variable at a time while keeping the others fixed. To support this iterative process and prove convergence, we define the iteration structure and introduce auxiliary variables as follows.
\begin{lstlisting}
class ADMM extends (Opt_problem E₁ E₂ F) where
  x₁ : ℕ → E₁   x₂ : ℕ → E₂  y  : ℕ → F
  ρ  : ℝ        τ  : ℝ       hrho : ρ > 0
  htau : 0 < τ ∧ τ < (1 + sqrt 5) / 2
  iterx₁  : ∀ k, IsMinOn (fun x₁ ↦ (Augmented_Lagrangian_Function E₁ E₂ F toOpt_problem ρ) (x₁, x₂ k, y k)) univ (x₁ (k + 1))
  uiterx₁ : ∀ k, Admm_sub_Isunique (fun x₁ ↦ (Augmented_Lagrangian_Function E₁ E₂ F toOpt_problem ρ) (x₁, x₂ k, y k)) (x₁ (k + 1)) (iterx₁ k)
  iterx₂  : ∀ k, IsMinOn (fun x₂ ↦ (Augmented_Lagrangian_Function E₁ E₂ F toOpt_problem ρ) (x₁ (k + 1), x₂, y k)) univ (x₂ (k + 1))
  uiterx₂ : ∀ k, Admm_sub_Isunique (fun x₂ ↦ (Augmented_Lagrangian_Function E₁ E₂ F toOpt_problem ρ) (x₁ (k + 1), x₂, y k)) (x₂ (k + 1)) (iterx₂ k)
  itery   : ∀ k, y (k + 1) = y k + (τ * ρ) • ((A₁ $ x₁ (k + 1)) + (A₂ $ x₂ (k + 1)) - b)
\end{lstlisting}

The parameters $\rho$ and $\tau$ range within the intervals $(0, \infty)$ and $\left(0, \frac{1 + \sqrt{5}}{2}\right)$, respectively. The property \lean{iterx₁} ensures that for any iteration index $k$, the point $x_1^{(k+1)}$ is a global minimizer of the function $L(x_1, x_2^{(k)}, y^{(k)})$. The uniqueness of this minimizer is established by \lean{uiterx₁}. Similarly, \lean{iterx₂} ensures that $x_2^{(k+1)}$ achieves global minimization for the function $L(x_1^{(k+1)}, x_2, y^{(k)})$, with its uniqueness established by \lean{uiterx₂}. The iteration of the dual variable $y$ is defined by \lean{itery}, which specifies the update formula for $y^{(k)}$ at each iteration.

The Karush-Kuhn-Tucker (KKT) conditions are necessary for optimality, given certain regularity assumptions. They are essential for determining the termination of the iterative process in the ADMM algorithm. For problem \eqref{eq: admm form}, these conditions have been adapted to align with the algorithm's framework and are expressed as follows:
\begin{equation}\label{eq : ADMM KKT}
  \begin{aligned}
  \begin{cases}
      &-A_1^Ty \in \nabla f_1(x_1), \\
     &-A_2^Ty \in \nabla f_2(x_2), \\
    &A_1x_1 + A_2x_2 = b.
  \end{cases}
  \end{aligned}
  \end{equation}
The formalization of the KKT conditions adheres to the structure of natural language and is presented below using the type \lean{class}.
\begin{lstlisting}
class Convex_KKT (x₁ : E₁) (x₂ : E₂) (y : F) (opt : Opt_problem E₁ E₂ F) : Prop where
  subgrad₁ : -(ContinuousLinearMap.adjoint opt.A₁) y ∈ SubderivAt opt.f₁ x₁
  subgrad₂ : -(ContinuousLinearMap.adjoint opt.A₂) y ∈ SubderivAt opt.f₂ x₂
  eq : (opt.A₁ x₁) + (opt.A₂ x₂) = opt.b
\end{lstlisting}

\subsection{Convergence Analysis}

This subsection describes the convergence of ADMM, with the convergence proof primarily based on \cite{Fazel2013hankelmatrix}. The proof is divided into two components: a lemma that verifies the descending property of an auxiliary sequence, and a theorem that establishes the convergence of the ADMM iterative sequence. We denote $(x_1^*, x_2^*, y^*)$ as the KKT pair for the original problem and define the error vectors and auxiliary functions as follows:
\[
\begin{aligned}
    (e_1^k,e_2^k,e_y^k) &= (x_1^k,x_2^k,y^k) - (x_1^*, x_2^*, y^*),\\
    u^k &= -A_1^T[y^k + (1-\tau)\rho(A_1e_1^k + A_2e_2^k) + \rho A_2(x_2^{k-1} - x_2^k)],\\
    v^k &= -A_2^T[y^k + (1-\tau)\rho(A_1e_1^k + A_2e_2^k)],\\
    \Psi_k &= \frac{1}{\tau\rho}\|e_y^k\|^2 + \rho\|A_2e_2^k\|^2,\\
    \Phi_k &= \Psi_k + \max(1-\tau, 1-\tau^{-1})\rho\|A_1e_1^k + A_2e_2^k\|^2.
\end{aligned}
\]
The functions $u^k$ and $v^k$ are related to the optimality conditions of each subproblem, $\Psi_k$ and $\Phi_k$ measure the errors $(e_1^k, e_2^k, y^k)$. All of these auxiliary variables are formalized using the \lean{ADMM} namespace. 

Our first step is to prove the descent of $\Phi_k$, implying the limits of the errors $(e_1^k, e_2^k, y^k)$ are zero.
\begin{lemma}\label{lemma: Lemma 20}
Suppose that Assumptions \ref{assumption: ADMM} and \ref{assumption: ADMM 2} hold. Let $\{(x_1^k,x_2^k,y^k)\}$ be a sequence generated by the ADMM iterations \eqref{eq: admm update 1} - \eqref{eq: admm update 3}. For each $k \geqslant 1$, it holds
\[
u^k \in \partial f_1(x_1^k),\quad v^k \in \partial f_2(x_2^k).
\]
The error obtains a decrease as 
\begin{align}\label{eq: admm theorem 1}
    \begin{aligned}
        \Phi_k - \Phi_{k+1} \geqslant & \min(\tau, 1+\tau-\tau^2)\rho\|A_2(x_2^k - x_2^{k+1})\|^2 \\&+ \min(1, 1+\tau^{-1}-\tau)\rho\|A_1e_1^{k+1} + A_2e_2^{k+1}\|^2.
    \end{aligned}
\end{align}

\end{lemma}
This lemma is formalized in Lean as below.
\begin{lstlisting}
lemma Φ_isdescending [Setting E₁ E₂ F admm admm_kkt] : ∀ n : ℕ+, (Φ n) - (Φ (n + 1)) ≥ (min τ (1 + τ - τ ^ 2)) * ρ * ‖A₂ (x₂ n - x₂ (n + 1))‖ ^ 2 + (min 1 (1 + 1 / τ - τ)) * ρ * ‖A₁ (e₁ (n + 1)) + A₂ (e₂ (n + 1))‖ ^ 2
\end{lstlisting}

We then outline the proof of this lemma, structured as follows:

\begin{enumerate}
  \item We prove that $u^k$ and $v^k$ belong to the subdifferentials of $f_1$ and $f_2$, respectively:
  \[
  \begin{aligned}
    u^{k+1} \in \partial f_1(x_1^{k+1}),  \quad  v^{k+1} \in \partial f_2(x_2^{k+1}).
  \end{aligned}
  \]
  This establishes that $u^k, v^k$ approximate the KKT pair and converge to it, due to the closed property of subgradients. The relevant proofs are provided in  \lean{uvBelongsToSubgradient}, \lean{u_inthesubgradient}, and \lean{v_inthesubgradient}.

  \item We establish an inequality that describes the lower bound of the difference in $\Phi_k$:
  \[
    \begin{aligned}
    &\frac{1}{\tau\rho} \left(\Vert e_y^k\Vert^2 - \Vert e_y^{k+1}\Vert^2\right)
    -(2 - \tau)\rho \Vert A_1x_1^{k+1} + A_2x_2^{k+1} - b\Vert^2\\
    &+2 M_{k+1} - \rho \Vert A_2(x_2^{k+1} - x_2^k)\Vert^2
    -\rho \Vert A_2e_2^{k+1}\Vert^2 + \rho \Vert A_2e_2^k\Vert ^2 \geqslant 0.
    \end{aligned}
  \]

  \item For cases where $0<\tau<1$ and $\tau > 1$, we categorize and discuss their properties. Finally, we can get the desired inequality \eqref{eq: admm theorem 1}.
\end{enumerate}
This lemma provides the foundation for proving convergence. It establishes that the iterations generated by ADMM remain bounded throughout the process.

With the lemmas established above, we need only demonstrate the convergence of the series $\{(x_1^k, x_2^k, y^k)\}$. This can be achieved by proving the convergence of subsequences and using the monotonicity of the sequence.
\begin{theorem}
Suppose that Assumptions \ref{assumption: ADMM} and \ref{assumption: ADMM 2} hold. Assuming that $A_1$ and $A_2$ have full rank and $\tau \in \left(0,\frac{1+\sqrt{5}}{2}\right)$, the series $\{(x_1^k, x_2^k, y^k)\}$ generated by the ADMM iterations \eqref{eq: admm update 1} - \eqref{eq: admm update 3} converge to a KKT pair of \eqref{eq: admm form}.
\end{theorem}
This theorem is defined in Lean as below.
\begin{lstlisting}
theorem ADMM_convergence [Setting E₁ E₂ F admm (admm_kkt₁ fullrank₁ fullrank₂ (admm_kkt := admm_kkt)] : ∃ (_x₁ : E₁) (_x₂ : E₂) (_y : F), Convex_KKT _x₁ _x₂ _y admm.toOpt_problem ∧ (Tendsto x₁ atTop (nhds _x₁) ∧ Tendsto x₂ atTop (nhds _x₂) ∧ Tendsto y atTop (nhds _y))
\end{lstlisting}

We hereby outline the key aspects of the proof, structured as follows:

\begin{enumerate}
  \item By Lemma \ref{lemma: Lemma 20}, the sequence $\{(x_1^k, x_2^k, y^k)\}$ is bounded. We can demonstrate that each of the sequences $x_1^k, x_2^k, y^k, u^k, v^k$ admits a convergent subsequence satisfying:
  \[
    \left(x_1^{\phi(k)}, x_2^{\phi(k)}, y^{\phi(k)}, u^{\phi(k)}, v^{\phi(k)}\right) \to (x_1'',x_2'',y'', -A_1^Ty'',-A_2^Ty''), \; \text{as} \; k \to \infty.
  \]
  Using monotonicity, we then prove the convergence of the entire sequence.

  \item The KKT conditions for problem \eqref{eq: admm form} are given by \eqref{eq : ADMM KKT}. Based on the boundedness of the iterates and the diminishing residuals, we can conclude that the sequence $\{x_1^k, x_2^k, y^k\}$ converges to a point satisfying the KKT conditions:
  \[
  \begin{aligned}
  \begin{cases}
    &-A_1^Ty'' \in \partial f_1(x_1''),\\
    &-A_2^Ty'' \in \partial f_2(x_2''),\\
    &A_1(x_1'') + A_2(x_2'') = b.
  \end{cases}
  \end{aligned}
  \]

  \item We derive the convergence of the entire sequences $x_1^k, x_2^k, y^k, u^k, v^k$ from the descending property of $\Phi^k$.
\end{enumerate}
Thus, under specified assumptions, the ADMM algorithm is guaranteed to converge to an optimal solution.

\section{Conclusion}
In this paper, the formalization of two first-order optimization algorithms, BCD and ADMM, specifically targeting block-structured optimization problems, are provided in Lean4. To address nonsmooth and nonconvex objectives, the subdifferential and the KL property are formalized, together with related theorems crucial for convergence analysis. We give the formal update scheme of the numerical algorithms using the structure type. Building on these concepts, we prove that the BCD algorithm converges to the critical points, and the ADMM algorithm converges to the KKT points of the problem. 

Our formalization framework can be naturally extended to derive further convergence results under alternative assumptions on the problem. Furthermore, numerous variations of BCD and ADMM documented in the literature can be formalized using our established methodology. For instance, various variants of the BCD algorithm, distinguished by slight modifications in their update schemes, can be systematically analyzed within this framework. Similarly, multiple versions of the ADMM algorithm are designed to accommodate specific problem structures. The formalization of these algorithms enables a more comprehensive study of the convergence properties of the ADMM method.




\bibliography{algorithm-formalization}

\begin{thebibliography}{10}

\bibitem{Appel2024VCFloat}
Andrew Appel and Ariel Kellison.
\newblock {VCFloat2: Floating-Point Error Analysis in Coq}.
\newblock In {\em Proceedings of the 13th ACM SIGPLAN International Conference on Certified Programs and Proofs}, CPP 2024, page 14–29, New York, NY, USA, 2024. Association for Computing Machinery.

\bibitem{attouch2010proximal}
H{\'e}dy Attouch, J{\'e}r{\^o}me Bolte, Patrick Redont, and Antoine Soubeyran.
\newblock {Proximal alternating minimization and projection methods for nonconvex problems: An approach based on the Kurdyka-{\L}ojasiewicz inequality}.
\newblock {\em Mathematics of operations research}, 35(2):438--457, 2010.

\bibitem{attouch2013convergence}
Hedy Attouch, J{\'e}r{\^o}me Bolte, and Benar~Fux Svaiter.
\newblock Convergence of descent methods for semi-algebraic and tame problems: proximal algorithms, forward--backward splitting, and regularized {G}auss--{S}eidel methods.
\newblock {\em Mathematical programming}, 137(1):91--129, 2013.

\bibitem{bento2024convergence}
GC~Bento, BS~Mordukhovich, TS~Mota, and Yu~Nesterov.
\newblock Convergence of descent methods under kurdyka-{\l}ojasiewicz properties.
\newblock {\em arXiv preprint arXiv:2407.00812}, 2024.

\bibitem{bolte2014proximal}
J{\'e}r{\^o}me Bolte, Shoham Sabach, and Marc Teboulle.
\newblock Proximal alternating linearized minimization for nonconvex and nonsmooth problems.
\newblock {\em Mathematical Programming}, 146(1-2):459--494, 2014.

\bibitem{boyd2011distributed}
Stephen Boyd, Neal Parikh, Eric Chu, Borja Peleato, Jonathan Eckstein, et~al.
\newblock Distributed optimization and statistical learning via the alternating direction method of multipliers.
\newblock {\em Foundations and Trends{\textregistered} in Machine learning}, 3(1):1--122, 2011.

\bibitem{chen2017note}
Liang Chen, Defeng Sun, and Kim-Chuan Toh.
\newblock A note on the convergence of admm for linearly constrained convex optimization problems.
\newblock {\em Computational Optimization and Applications}, 66:327--343, 2017.

\bibitem{de2015lean}
Leonardo De~Moura, Soonho Kong, Jeremy Avigad, Floris Van~Doorn, and Jakob von Raumer.
\newblock The lean theorem prover (system description).
\newblock In {\em Automated Deduction-CADE-25: 25th International Conference on Automated Deduction, Berlin, Germany, August 1-7, 2015, Proceedings 25}, pages 378--388. Springer, 2015.

\bibitem{dupuis_et_al:LIPIcs.ITP.2022.10}
Fr\'{e}d\'{e}ric Dupuis, Robert~Y. Lewis, and Heather Macbeth.
\newblock {Formalized functional analysis with semilinear maps}.
\newblock In June Andronick and Leonardo de~Moura, editors, {\em 13th International Conference on Interactive Theorem Proving (ITP 2022)}, volume 237 of {\em Leibniz International Proceedings in Informatics (LIPIcs)}, pages 10:1--10:19, Dagstuhl, Germany, 2022. Schloss Dagstuhl -- Leibniz-Zentrum f{\"u}r Informatik.

\bibitem{Fazel2013hankelmatrix}
Maryam Fazel, Ting~Kei Pong, Defeng Sun, and Paul Tseng.
\newblock Hankel matrix rank minimization with applications to system identification and realization.
\newblock {\em SIAM Journal on Matrix Analysis and Applications}, 34(3):946--977, 2013.

\bibitem{fortin2000augmented}
Michel Fortin and Roland Glowinski.
\newblock {\em Augmented Lagrangian methods: applications to the numerical solution of boundary-value problems}, volume~15.
\newblock Elsevier, 2000.

\bibitem{gabay1976dual}
Daniel Gabay and Bertrand Mercier.
\newblock A dual algorithm for the solution of nonlinear variational problems via finite element approximation.
\newblock {\em Computers \& mathematics with applications}, 2(1):17--40, 1976.

\bibitem{glowinski2013numerical}
Roland Glowinski.
\newblock {\em Numerical methods for nonlinear variational problems}.
\newblock Springer Science \& Business Media, 2013.

\bibitem{Glowinski2014}
Roland Glowinski.
\newblock {\em On Alternating Direction Methods of Multipliers: A Historical Perspective}, pages 59--82.
\newblock Springer Netherlands, Dordrecht, 2014.

\bibitem{glowinski1975approximation}
Roland Glowinski and Americo Marroco.
\newblock Sur l'approximation, par {\'e}l{\'e}ments finis d'ordre un, et la r{\'e}solution, par p{\'e}nalisation-dualit{\'e} d'une classe de probl{\`e}mes de dirichlet non lin{\'e}aires.
\newblock {\em Revue fran{\c{c}}aise d'automatique, informatique, recherche op{\'e}rationnelle. Analyse num{\'e}rique}, 9(R2):41--76, 1975.

\bibitem{huet1997coq}
G{\'e}rard Huet, Gilles Kahn, and Christine Paulin-Mohring.
\newblock The {Coq} proof assistant a tutorial.
\newblock {\em Rapport Technique}, 178, 1997.

\bibitem{Kellison2024Numerical}
Ariel~E. Kellison and Justin Hsu.
\newblock {Numerical Fuzz: A Type System for Rounding Error Analysis}.
\newblock {\em Proc. ACM Program. Lang.}, 8(PLDI), June 2024.

\bibitem{li2024formalization}
Chenyi Li, Ziyu Wang, Wanyi He, Yuxuan Wu, Shengyang Xu, and Zaiwen Wen.
\newblock {Formalization of Complexity Analysis of the First-order Algorithms for Convex Optimization}, 2024.
\newblock \href {https://arxiv.org/abs/2403.11437} {\path{arXiv:2403.11437}}.

\bibitem{lojasiewicz1963propriete}
Stanislaw Lojasiewicz.
\newblock Une propri{\'e}t{\'e} topologique des sous-ensembles analytiques r{\'e}els.
\newblock {\em Les {\'e}quations aux d{\'e}riv{\'e}es partielles}, 117:87--89, 1963.

\bibitem{mathlibcommunity}
The mathlib Community.
\newblock The {L}ean mathematical library.
\newblock In {\em Proceedings of the 9th ACM SIGPLAN International Conference on Certified Programs and Proofs}, CPP 2020, page 367–381, New York, NY, USA, 2020. Association for Computing Machinery.

\bibitem{Nipkow2002APA}
Tobias Nipkow, Lawrence~Charles Paulson, and Markus Wenzel.
\newblock A proof assistant for higher-order logic.
\newblock {\em Lecture Notes in Computer Science}, 2002.

\bibitem{polyak1963gradient}
Boris~T Polyak.
\newblock Gradient methods for the minimisation of functionals.
\newblock {\em USSR Computational Mathematics and Mathematical Physics}, 3(4):864--878, 1963.

\bibitem{powell1973search}
M.J.D. Powell.
\newblock On search directions for minimization algorithms.
\newblock {\em Mathematical Programming}, 4(1):193--201, 1973.

\bibitem{Tekriwal2023Verified}
Mohit Tekriwal, Andrew~W. Appel, Ariel~E. Kellison, David Bindel, and Jean-Baptiste Jeannin.
\newblock {Verified Correctness, Accuracy, and Convergence of a Stationary Iterative Linear Solver: Jacobi Method}.
\newblock In {\em Intelligent Computer Mathematics: 16th International Conference, CICM 2023, Cambridge, UK, , September 5–8, 2023 Proceedings}, page 206–221, Berlin, Heidelberg, 2023. Springer-Verlag.

\bibitem{Tekriwal2024Formalizatioin}
Mohit Tekriwal, Joshua Miller, and Jean-Baptiste Jeannin.
\newblock {Formalization of Asymptotic Convergence for Stationary Iterative Methods}.
\newblock In {\em NASA Formal Methods: 16th International Symposium, NFM 2024, Moffett Field, CA, USA, June 4–6, 2024, Proceedings}, page 37–56, Berlin, Heidelberg, 2024. Springer-Verlag.

\bibitem{tseng2001convergence}
Paul Tseng.
\newblock Convergence of a block coordinate descent method for nondifferentiable minimization.
\newblock {\em Journal of optimization theory and applications}, 109:475--494, 2001.

\bibitem{vandoorn_et_al:LIPIcs.ITP.2024.37}
Floris van Doorn and Heather Macbeth.
\newblock {Integrals Within Integrals: A Formalization of the Gagliardo-Nirenberg-Sobolev Inequality}.
\newblock In Yves Bertot, Temur Kutsia, and Michael Norrish, editors, {\em 15th International Conference on Interactive Theorem Proving (ITP 2024)}, volume 309 of {\em Leibniz International Proceedings in Informatics (LIPIcs)}, pages 37:1--37:18, Dagstuhl, Germany, 2024. Schloss Dagstuhl -- Leibniz-Zentrum f{\"u}r Informatik.

\bibitem{Ying2023Doob}
Kexing Ying and R\'{e}my Degenne.
\newblock {A Formalization of Doob’s Martingale Convergence Theorems in mathlib}.
\newblock In {\em Proceedings of the 12th ACM SIGPLAN International Conference on Certified Programs and Proofs}, CPP 2023, page 334–347, New York, NY, USA, 2023. Association for Computing Machinery.

\end{thebibliography}

\end{document}